\documentclass{amsart}
\usepackage{amsfonts}
\usepackage{amsmath}
\usepackage{amssymb}
\usepackage{amscd}
\usepackage{amstext}
\usepackage{amsthm}
\usepackage{enumerate}
\usepackage{color}

\theoremstyle{plain}
\newtheorem{lem}{Lemma}[section]

\newtheorem{thm}[lem]{Theorem}

\theoremstyle{definition}

\theoremstyle{remark}

\begin{document}
\title{Duistermaat-Heckman formula for a torus action on a generalized Calabi-Yau manifold and localization formula}
\author{Yasufumi Nitta}
\address{Department of Mathematics, Graduate school of science, Osaka University. 1-1 Machikaneyama, Toyonaka, Osaka 560-0043 Japan.}

\maketitle

\begin{abstract}
This note is an addendum to our earlier work \cite{humi}. In \cite{humi}, we studied a Hamiltonian action for 
a generalized Calabi-Yau manifold and showed that the Duistermaat-Heckman theorem holds. The purpose of this note is 
to show that the density function of the Duistermaa-Heckman measure is a piecewise polynomial. We also prove that 
the localization formula holds.  
\end{abstract}

\section{introduction}
In \cite{humi}, we studied a Hamiltonian action for a generalized Calabi-Yau manifold. Let $M$ be a $2n$-dimensional connected 
generalized Calabi-Yau manifold with generalized Calabi-Yau structure $\varphi$ of constant type $k$, 
and let $T$ be a $l$-dimensional compact torus which acts on $M$ in a Hamiltonian way. 
That is, this is a Lie group action on $M$ which preserves the generalized Calabi-Yau structure $\varphi$, and there is a smooth function 
$\mu : M \longrightarrow \mathfrak{t}^{*}$ which satisfies 
\begin{enumerate}[{\rm (1)}]
\item $\mu$ is $T$-equivariant, and 
\item $\xi_{M} - \sqrt{-1}d\mu^{\xi}$ lies in $E_{\varphi}$ for all $\xi \in \frak{t}$, where $\xi_{M}$ denotes the induced vector field on $M$ and 
$\mu^{\xi}$ is the smooth function defined by $\mu^{\xi}(p) = \mu (p)(\xi)$.
\end{enumerate}
We call the function $\mu$ the generalized moment map for the Hamiltonian action. Let $dm$ be a natural volume form on $M$ defined by 
\[
dm = \frac{(\sqrt{-1})^n}{2^{n-k}} \langle \varphi , \bar \varphi \rangle. 
\]
We proved in \cite{humi} that the pushforward measure $\mu_{*}(dm)$ on $\mathfrak{t}^{*}$ is absolutely continuous with respect to 
the Lebesgue measure $d\zeta$ on $\mathfrak{t}^{*}$, and gave the explicit formula of the density function. 

\begin{thm}[\cite{humi}, Theorem B.]
Let $(M, \varphi)$ be a $2n$-dimensional connected generalized Calabi-Yau manifold of constant type k, and suppose that compact $l$-torus $T$ 
acts on $M$ effectively and in a Hamiltonian way. In addition, We assume that the generalized moment map $\mu$ is proper. 
Then the pushforward $\mu_{*}(dm)$ of the natural volume form $dm$ under $\mu$ is absolutely continuous with respect to the 
Lebesgue measure on $\frak{t}^{*}$ and the Radon-Nikodym derivative $f$ can be written by 
\[
f(a) = \int_{M_{a}}dm_{a} = {\rm vol}(M_{a}) 
\]
for each regular value $a \in \frak{t}^{*}$ of $\mu$, where 
$dm_{a}$ denotes the measure defined by the natural volume form on the reduced space $M_{a} = \mu^{-1}(a)/T$. 
\end{thm}

The main purpose of this note is to prove the density function $f$ is a piecewise polynomial. 

\begin{thm}\label{main1}
Under the assumptions above, 
the function $f$ is a polynomial (of degree $\leq  n-k-l$) on each connected component of the set of regular values of 
the generalized moment map $\mu$. 
\end{thm}
The proof of Theorem \ref{main1} will be given in Section 2. Next using the Atiyah-Bott-Berline-Vergne localization theorem, 
we shall lead the localization formula for the oscillatory integral
\[
\int_{M}e^{\sqrt{-1}\mu^{\xi}}dm,\quad \xi \in \frak{t}
\]
in the case that $M$ is compact. 
The detailed notations and definitions is in the Section 3. 

\begin{thm}\label{main2}
In addition to the assumptions in Theorem \ref{main1}, we assume that $M$ is compact. Then the localization theorem holds, i.e., 
the oscillatory integral is given by 
\[
\int_{M}e^{\sqrt{-1}\mu^{\xi}}dm = \sum_{X \subset M^{T}}\frac{1}{2^{(\dim X - l)}}\int_{X}\frac{e^{i_{X}^{*}\mu^{\xi}}i_{X}^{*}(\varphi \wedge \sigma ( \bar \varphi ))}{e(NX)(\xi)}, 
\]
for $\xi \in \frak{t}$ such that $\beta_{j}(\xi) \not = 0$, where $X$ is a connected component of the fixed point set $M^{T}$, $i_{X} : X \longrightarrow M$ is the inclusion, and 
$e(NX)$ is the equivariant Euler class of the normal bundle of $X$. 
\end{thm}

\section{Proof of Theorem \ref{main1}}
Let $M_{{\rm free}}$ be the set on which $T$ acts freely. 
Then since $M$ is connected and $T$ acts on $M$ effectively, $M_{{\rm free}}$ is open dense subset of $M$ and $M \setminus M_{{\rm free}}$ 
has measure $0$ (see \cite{Gui_Stern}, Appendix B). 
So we may assume $T$ acts $M$ freely, and it is sufficient to show that $f$ is a polynomial of degree at most $n-l-k$. 
Let $p : M \longrightarrow M/T$ be the natural projection. Since the generalized moment map $\mu : M \longrightarrow \frak{t}^{*}$ 
is $T$-invariant, there exists a smooth map 
\[
\psi : M/G \longrightarrow \frak{t}^{*}
\]
such that 
\[
\mu = \psi \circ p. 
\]
Then $\mu^{-1}(a) = p^{-1}(\psi^{-1}(a))$ and $M_{a} = \psi^{-1}(a)$, for each regular value $a \in \mathfrak{t}^{*}$. 
Let $j_{a} : M_{a} \longrightarrow M/G$ be the inclusion, so that we have the commutative relation
\[
\pi \circ i_{a} = j_{a} \circ p_{a}, 
\]
where $i_{a} : \mu^{-1}(a) \longrightarrow M$ is the inclusion, and $p_{a} : \mu^{-1}(a) \longrightarrow M_{a}$ is the natural projection. 
Now we choose a connection on $M$ and let $\theta^{i}$ be the curvature forms. Then the equivariant differential form
\[
e^{2\sqrt{-1}\mu}\varphi \wedge \sigma(\bar \varphi)
\]
is closed, and using the Cartan map $\Omega_{T}(M) \longrightarrow \Omega(M/T)$, 
\[
e^{2\sqrt{-1}\psi_{i}\theta^{i}} \wedge \nu \wedge \sigma(\bar \nu) 
\]
is closed form, where $\nu$ is the unique form on $M/T$ with the property that 
\[
p^{*}\nu = \varphi_{{\rm hor}}, 
\]
the horizontal part of $\varphi$ (the definition of the Cartan map is in \cite{Gui_Stern}, for instance). From the relation $\pi \circ i_{a} = j_{a} \circ \pi_{a}$, we obtain 
\[
j_{a}^{*}\nu = \varphi_{a}. 
\]
Let $c$ denote the cohomology class of $e^{2\sqrt{-1}\psi_{i}\theta^{i}}\nu \wedge \sigma(\bar \nu)$, then 
\begin{eqnarray*}
[\varphi_{a} \wedge \sigma(\bar\varphi_{a})] &=& j_{a}^{*}[\nu \wedge \sigma(\bar \nu) ] \\
                                             &=& j_{a}^{*}[ e^{-2\sqrt{-1}\psi_{i}\theta^{i}} \wedge (e^{2\sqrt{-1}\psi_{i}\theta^{i}} \wedge \nu \wedge \sigma(\bar \nu) ) ] \\
                                             &=& j_{a}^{*} \left( e^{-2\sqrt{-1}a_{i}[\theta^{i}]}\wedge c \right). 
\end{eqnarray*}
Since $M_{a}$ is compact and oriented, the embedding $j_{a} : M_{a} \longrightarrow M/T$ defines a homology class $[M_{a}] \in H_{2(n - l)}(M/T, \mathbb{Z})$. This depends smoothly 
on $a$, and being an integer class is thus independent of $a$. So let us fix an $a_{0}$ in the image of $\mu$. Then we have $[M_{a_{0}}] = [M_{a}]$. 
Thus the integral can be calculated by 
\begin{eqnarray*}
f(a) &=& \left\langle [M_{a}],\ [dm_{a}] \right\rangle \\
     &=& \left\langle [M_{a_{0}}],\ [\frac{(\sqrt{-1})^{(n-l)}}{2^{n-l-k}} (\varphi_{a} \wedge \sigma (\bar \varphi_{a}))_{2(n-l)} ] \right\rangle \\
     &=& \left\langle [M_{a_{0}}],\ \frac{(\sqrt{-1})^{(n-l)}}{2^{n-l-k}}( e^{-2\sqrt{-1}a_{i}[\theta^{i}]}\wedge c )_{2(n-l)} \right\rangle, 
\end{eqnarray*}
where $\langle , \rangle$ is the pairing of a homology class and a cohomology class. 
This shows that $f$ is a polynomial of degree at most $n-l-k$. 

\section{Proof of Theorem \ref{main2}}
Suppose that $M$ is compact. Then the number of subgroups of $T$ which can occur as 
stabilizer groups of points of $M$ is finite, and so there exists a finite number of weights 
$\beta_{1}, \cdots , \beta_{m}$ with the property that if $\xi \in \mathfrak{t}$ satisfies 
\[
\beta_{i}(\xi) \not = 0,\ i = 1, \cdots , m
\]
the corresponding vector field $\xi_{M}$ on $M$ is non-zero except at points of $M^{T}$(see \cite{Gui_Stern} for example). 

Consider a complex form $\varphi^{\prime} \in \Omega^{ev/od}M \otimes \mathbb{C}$ defined by 
\[
\varphi^{\prime} = \varphi_{k} + \frac{1}{2}\varphi_{k+2} + \frac{1}{4}\varphi_{k+4} + \cdots, 
\]
where $\varphi_{s}$ is the $s$-th degree component of $\varphi$. Then a complex equivariant differential form 
$e^{\sqrt{-1}\mu}\varphi^{\prime} \wedge \sigma(\bar{\varphi^{\prime}})$ is equivariant closed, and 
$(\varphi^{\prime} \wedge \sigma(\bar{\varphi^{\prime}}))_{2n} = dm$. So by applying the Atiyah-Bott-Berline-Vergne localization theorem
(see for instance \cite{At} and \cite{Gui_Stern}),

we have 
\begin{eqnarray*}
\int_{M}e^{\sqrt{-1}\mu^{\xi}}dm &=& \int_{M}(e^{\sqrt{-1}\mu}\varphi^{\prime} \wedge \sigma(\bar{\varphi^{\prime}}))(\xi) \\
                                                  &=& \sum_{X \subset M^{T}}\int_{X}\frac{i_{X}^{*}(e^{\sqrt{-1}\mu}\varphi^{\prime} \wedge \sigma(\bar{\varphi^{\prime}})))(\xi)}{e(NX)(\xi)} \\
                                                  &=& \sum_{X \subset M^{T}}\frac{1}{2^{(\dim X - l)}}\int_{X}\frac{e^{i_{X}^{*}\mu^{\xi}}i_{X}^{*}(\varphi \wedge \sigma ( \bar \varphi ))}{e(NX)(\xi)},  \\
\end{eqnarray*}
for $\xi \in \frak{t}$ satisfying $\beta_{i}(\xi) \not = 0$ for all $i = 1, \cdots , m$, where $X$ is a connected component of the fixed point set $M^{T}$, 
$i_{X} : X \longrightarrow M$ is the inclusion, and $e(NX)$ is the equivariant Euler class of the normal bundle of $X$. This shows Theorem \ref{main2}.

\end{document}